\documentclass[graybox]{svmult}
\usepackage{type1cm}        
\usepackage{makeidx}         
\usepackage{graphicx}        
\usepackage{multicol}        
\usepackage[bottom]{footmisc}
\usepackage{multicol,url}
\usepackage{latexsym, amsmath, amssymb}
\usepackage{amsfonts}
\usepackage{pifont}
\usepackage{enumitem}
\usepackage[colorlinks=true]{hyperref}
\usepackage{empheq}
\usepackage{tabularx,multirow,array,subcaption,float} 
\usepackage{xparse}
\ExplSyntaxOn
\int_new:N \g_tohi_const_int
\int_new:N \g_tohi_const_sub_int
\tl_new:N  \g_tohi_const_char_tl

\cs_new_protected:Nn \tohi_print_constant:nn 
{
  #1 \textsubscript {#2} 
} 

\NewDocumentCommand\resetconstants{m}
{
 \int_gincr:N \g_tohi_const_int
 \int_gzero:N \g_tohi_const_sub_int
 \tl_gset:Nn  \g_tohi_const_char_tl {#1}
}

\NewDocumentCommand\const{m}
{
  \tl_if_exist:cTF
   {
    c_tohi_const_\int_use:N\g_tohi_const_int _#1_tl
   }
   {
    \tl_use:c {c_tohi_const_\int_use:N\g_tohi_const_int _#1_tl }
   }
   {
    \int_gincr:N \g_tohi_const_sub_int
    \tl_const:cx {c_tohi_const_\int_use:N\g_tohi_const_int _#1_tl }
     { \exp_not:N\tohi_print_constant:nn {\g_tohi_const_char_tl }{\int_use:N \g_tohi_const_sub_int}}
    \tl_use:c {c_tohi_const_\int_use:N\g_tohi_const_int _#1_tl }
   }
}

\ExplSyntaxOff
\newcommand{\inlineitem}[1][]{%
\ifnum\enit@type=\tw@
    {\descriptionlabel{#1}}
  \hspace{\labelsep}%
\else
  \ifnum\enit@type=\z@
       \refstepcounter{\@listctr}\fi
    \quad\@itemlabel\hspace{\labelsep}%
\fi}
\DeclarePairedDelimiter\abs{\lvert}{\rvert}
\DeclarePairedDelimiter\norm{\lVert}{\rVert}
\DeclarePairedDelimiter\tonda{(}{)}

\DeclarePairedDelimiter\graffa{\{}{\}}
\newcommand{\into}{\int_\Omega}

\setlist[itemize]{noitemsep, topsep=0pt}
\def\R{\mathbb R} \def\N{\mathbb N} 
\def\tmax{\ensuremath{T_{\textup{max}}}}

\def\R{\mathbb R} \def\N{\mathbb N} 
\def\TM{T_{\textup{max}}}

\def\dqvq{\abs*{D^2 v}^2}
\def\absnv{\abs{\nabla v}}
\def\nivz{\norm{v_0}_{L^{\infty}(\Omega)}}

\def
\@cite
#1#2{[#1\if@tempswa, #2\fi]}
\makeatother

\makeindex 
\usepackage{newtxtext}       %
\usepackage[varvw]{newtxmath}       
\begin{document}
\title*{Boundedness in a nonlinear chemotaxis-consumption model with gradient terms}
\author{Daniel Acosta-Soba\orcidID{0000-0002-2159-5926} \\ Alessandro Columbu\orcidID{0000-0001-6993-1223}  and\\ Giuseppe Viglialoro\orcidID{0000-0002-2994-4123}}
\institute{Daniel Acosta-Soba \at Departamento de Matemáticas, Universidad de Cádiz, Puerto Real, Campus Universitario Río San Pedro s/n, 11510.  Cádiz (Spain)
\email{daniel.acosta@uca.es}
\and Alessandro Columbu (\emph{Corresponding author})\at Dipartimento di Matematica e Informatica, Universit\`{a} di Cagliari, Via Ospedale 72, 09124.\\ Cagliari (Italy) \email{alessandro.columbu2@unica.it}
\and Giuseppe Viglialoro \at Dipartimento di Matematica e Informatica, Universit\`{a} di Cagliari, Via Ospedale 72, 09124.\\ Cagliari (Italy)  \email{giuseppe.viglialoro@unica.it}}
\maketitle
%
%
\abstract{We study a chemotaxis-consumption mechanism, in which some chemical signal and cells density interact each other. In order to control the concentration of such a population, sources  involving gradient nonlinearities, which introduce a dampening effect on the model, are considered. Moreover, the system is characterized by nonlinear diffusion and sensitivity terms. We derive conditions on some data of the problem so to ensure the boundedness of related solutions. This work extends the research presented in \cite{MarrasViglialoroMathNach,Columbu_arXiv_2024},  where the same nonlinear model without gradient terms and its linear version with gradient sources has been, respectively, addressed.}
\resetconstants{c}
\section{Introduction and motivations}\label{IntroductionSection}
\subsection{The landmarking Keller--Segel models with production and  consumption}
Systems of interacting agents are widespread in the physical and biological sciences. Examples range from predator-prey models to taxis-driven pattern formation and front propagation in mathematical biology. In this context, we refer to the second landmarking paper by Keller and Segel (\cite{Keller-1971-TBC}), which involves two coupled parabolic equations. The mathematical formulation of this system is given as the following initial-boundary value problem:
\begin{equation}\label{Keller-Segel-CosnumptionTao}
\left\{ 
\begin{array}{ll}
u_{t} =  \Delta u - \chi \nabla \cdot ( u  \nabla v) & \text{in } \Omega \times (0, \TM), \\
v_{t} = \Delta v - uv & \text{in } \Omega \times (0, \TM), \\
u_\nu = v_\nu = 0 & \text{on } \partial \Omega \times (0, \TM), \\
u(x, 0) = u_0(x) \geq 0 \quad \text{and} \quad v(x, 0) = v_0(x) \geq 0 & x \in \bar{\Omega}.
\end{array} 
\right. 
\end{equation}
Here, $\chi > 0$, $\Omega \subset \mathbb{R}^n$, with $n \geq 2$, is a bounded domain with smooth boundary $\partial \Omega$ and $\TM\in (0,\infty]$ the maximum instant of time up to which model evolves.  This problem describes the situation where the motion of certain individual cells, represented by $u = u(x,t)$, is influenced by the presence of a chemical signal, $v = v(x,t)$, which, due to the negative term $-uv$ in the second equation, is progressively consumed by the cells. Naturally, this suggests that $v$ remains bounded over time. 

Additionally, $u_0(x)$ and $v_0(x)$ represent the initial cell and chemical distributions, respectively. The symbol $(\cdot)_\nu$  denotes the outward normal derivative on $\partial \Omega$, and the zero-normal boundary conditions on both $u$ and $v$ imply that the domain is completely insulated.

The described phenomenon is quite different from the classical Keller--Segel model (see \cite{K-S-1970}), where the term $-uv$ is replaced by $-v + u$. The corresponding formulation reads (avoiding to insert intial and boundary conditions) as follows:
\begin{equation}\label{Keller-Segel-calssic} 
\left\{ 
\begin{array}{ll}
u_{t} = \Delta u - \chi \nabla \cdot ( u  \nabla v) & \text{in } \Omega \times (0, \TM), \\
v_{t} = \Delta v - v + u & \text{in } \Omega \times (0, \TM).
\end{array} 
\right. 
\end{equation}
This formulation highlights how an increase in the ammount of cells leads to the production of the chemical signal. Consequently, no a priori bound for $v$ is expected. Although deeply related, the two models \eqref{Keller-Segel-CosnumptionTao} and \eqref{Keller-Segel-calssic} exhibit different properties. In particular, for system \eqref{Keller-Segel-calssic}, a comprehensive and extensive theory exists concerning the existence and properties of global  (i.e., $\TM=\infty$) solutions, as well as uniformly bounded or blow-up solutions (which become unbounded in finite or infinite time $\TM$). A thorough treatment of this theory is available in the introduction of \cite{HorstWink}, and further surveys can be found in the work by Hillen and Painter \cite{Hillen-2009-UGP}, which includes reviews of modeling issues in various Keller--Segel-type systems.

On the other hand, for  system \eqref{Keller-Segel-CosnumptionTao}, it has been shown that the existence of global or blow-up solutions is independent of the size if initial data $u_0$. Specifically, Tao (\cite{TaoBoun}) proves that for sufficiently regular initial data $u_0$ and $v_0$, if 
\begin{equation}\label{Conditionv0chiTao}  
0 < \chi \lVert v_0 \rVert_{L^\infty(\Omega)} \leq \frac{1}{6(n+1)},
\end{equation}
then problem \eqref{Keller-Segel-CosnumptionTao} possesses a unique global classical solution that remains uniformly bounded. (See \cite{BaghaeiKhelghati-nolog} for an improvement of this condition.) The question of whether blow-up solutions exist for large initial data $v_0$ or large chemotactic parameter $\chi$, which do not satisfy \eqref{Conditionv0chiTao}, remains open.

Further progress on the chemotaxis-consumption model \eqref{Keller-Segel-CosnumptionTao} has been made. In \cite{TaoWinkConsumptionEventual}, weak solutions that become smooth after some time are constructed for sufficiently regular initial data in the three-dimensional setting. For $n = 2$, globally bounded classical solutions with convergence properties are also shown to exist.

Moreover, by interpreting the system as a special case of the general coupled chemotaxis-fluid model proposed by Goldstein in \cite{TuvalGoldsteinEtAl}, where the fluid does not contribute directly, \cite{WinklerN-Sto_CPDE} discusses the existence of global classical and weak solutions for $n = 2$ and $n = 3$, respectively. In addition, \cite{WinklerN-Sto_2d} investigates the stabilization properties of the two-dimensional solutions. Further results on existence and asymptotic behavior of solutions for the three-dimensional version of the coupled system are provided in \cite{WinklerGlobalSolNvStokes} and \cite{WinklerAMS}.
\subsection{The introduction of logistic degradation}
In order to contrast the undesired blow-up singularities, which, as mentioned above, may arise in both models \eqref{Keller-Segel-CosnumptionTao} and \eqref{Keller-Segel-calssic}, more complete formulations of these systems with nontrivial sources have been considered. Specifically, it is natural to complement these models with \textit{logistic-type} effects, which are widely used in mathematical biology. 

For instance, let us consider the system 
\begin{equation*}\label{sys Johannes}
\begin{cases}
u_{t} = \Delta u -  \chi \nabla \cdot (u \nabla v) + g(u) & \text{in } \Omega \times (0, \TM), \\
\tau v_{t} = \Delta v - v + u & \text{in } \Omega \times (0, \TM),
\end{cases}
\end{equation*}
 endowed with homogeneous Neumann boundary conditions and sufficiently regular initial data, and  perturbed by a logistic-type source term, $g(u) \simeq \lambda u - \mu u^\gamma$, with $\lambda  \in \mathbb{R}$ and $\mu > 0$ as constants. 

In \cite{Lankeit}, the existence of global weak solutions is established for the case $\chi = \tau = 1$ and $\gamma = 2$, for any arbitrarily small value of $\mu > 0$. Moreover, if $n = 3$, these solutions become classical after some time, provided that $k$ is not too large. On the other hand, for $\mu, \tau > 0$, $\gamma = 2$, and $\chi >0$, the author in \cite{W0} proves that if $\mu$ is sufficiently large, the same system possesses a unique bounded and global-in-time classical solution.

Finally, for $\gamma > 1$, $\tau = 1$, and $\chi > 0$, the global existence of very weak solutions, as well as their boundedness properties and long-time behavior, are discussed in \cite{ViglialoroVeryWeak}, \cite{ViglialoroBoundnessVeryWeak}, and \cite{ViglialoroWolleyDCDS}.

As explicitly derived in \cite{BelBelNieSol}, generalized kinetic theory methods can be intrinsically employed to obtain the continuous models presented here from the underlying description at the scale of cells. In this sense, the evolution of $u$, which has so far been considered (and idealizes a diffusion with infinite speed of propagation, which is definitely not lifelike), is a special case of the equation
\begin{equation*}
u_t = \nabla \cdot \left(H(u,v)\nabla u - K(u,v)\nabla v \right) \quad \text{in } \Omega \times (0, \TM).
\end{equation*}
This latter formulation includes a very general framework where the diffusion  $H$ and sensitivity $K$ may be described by different expressions. For instance, some special cases with singular sensitivity are noteworthy. Specifically, in the one-dimensional setting, if $H(u,v)$ is a positive constant and $K(u,v) = \frac{\chi u}{v}$, results on global existence for initial-boundary value problems tied to \eqref{Keller-Segel-CosnumptionTao} with arbitrary initial data have been derived in \cite{TaoWangEtALDCDS-B}. In higher dimensions, the analysis is more complex, and only partial results are known (see \cite{WinklerSingularSensitivityRelaxation} and \cite{win_ct_sing_abs_renormalized} for interesting achievements involving generalized and renormalized solutions in bounded domains of $\mathbb{R}^n$, $n \geq 2$).

Additionally, in \cite{LankeitLocallyBoundedSingularity} and \cite{dongmeiGlo}, with the respective choices $H(u,v) \geq \delta u^{m_1-1}$, $\delta > 0$, and $K(u,v) = \frac{u}{v}$, and $H(u,v) \equiv 1$ and $K(u,v) \simeq \frac{\lambda u (u+1)^{m_2-1}}{v}$, $\lambda > 0$, for appropriate values of $m$ and $m_2$ (depending on the dimension $n$), global classical solutions have been studied. On the other hand, for the problem with $H(u,v) = (u+1)^{m_1-1}$ and $K(u,v) = u(u+1)^{m_2-1}$, for some $m_1,m_2 \in \mathbb{R}$, and the equation for $v$ as in \eqref{Keller-Segel-calssic}, it has been studied in \cite{CieslakStinnerFiniteTime}, \cite{CieslakStinnerNewCritical}, and \cite{TaoWinkParaPara}. It is essentially established that the condition $m_2 < m_1 + \frac{2}{n} - 1$ is both sufficient and necessary to ensure global existence and boundedness of solutions, even for large initial data. 

To the best of our knowledge, significant results in the direction of this investigation, tied to \eqref{Keller-Segel-CosnumptionTao} under a perturbation for the evolution of $u$ through a logistic source, are as follows. In a bounded and smooth domain of $\mathbb{R}^n$, $n \geq 1$, and under homogeneous Neumann boundary conditions, for the system
\begin{equation}\label{sysGeneralForConsumptionAndLogisticJohannes}
\begin{cases}
u_t = \nabla \cdot (D(u) \nabla u) - \nabla \cdot (u \chi(u,v) \nabla v) + g(u) & \text{in } \Omega \times (0, \TM), \\
v_t = \Delta v - u f(v) & \text{in } \Omega \times (0, \TM),       
\end{cases}
\end{equation}
these results have been achieved (the list is not exhaustive):
\begin{enumerate}[label=\roman*)]
\item \label{BaghaeiaKhelghatibItems} For $D(u) = \delta > 0$, $\chi(u,v) = \chi > 0$, and $g(u) \leq \lambda u - \mu u^\gamma$ (with $g(0) \geq 0$, $\lambda \geq 0$, $\mu > 0$, and $\gamma > 1$), and sufficiently regular initial data $(u_0, v_0)$, the system admits a unique global and bounded classical solution for suitable small $\chi \lVert v_0 \rVert_{L^\infty(\Omega)}$ (see \cite{BaghaeiaKhelghatib}).
\item \label{LankeitConsumptionItem} For $D(u) \equiv 1$, $\chi(u,v) = \chi > 0$, $\lambda \in \mathbb{R}$, $g(u) = \lambda u - \mu u^2$, and sufficiently regular initial data $(u_0, v_0)$, the system admits a unique global and bounded classical solution for $\mu$ large enough compared to $\chi \lVert v_0 \rVert_{L^\infty(\Omega)}$, and a weak solution for arbitrary $\mu > 0$ (see \cite{LankeitWangConsumptLogisticabbr}).
\item \label{MarrasViglialoroItem} For $D(u) \equiv (u + 1)^{m_1-1}$, $\chi(u,v) = \chi (u + 1)^{m_2-1}$ ($\chi > 0$), $m_1, m_2, \lambda \in \mathbb{R}$, $g(u) = \lambda u - \mu u^2$, $f(v)=v$, it is proved that for nonnegative and sufficiently regular initial data \( u(x, 0) \) and \( v(x, 0) \), the corresponding initial-boundary value problem admits a unique globally bounded classical solution, provided $m_2<\frac{m_1+1}{2}$ and \( \mu \) is larger than a quantity depending on $\chi \lVert v_0\rVert_{L^\infty(\Omega)}$ (see \cite{MarrasViglialoroMathNach}).  
\end{enumerate}
\begin{remark}\label{RemarckChiV_0}
The contributions summarized in items \ref{BaghaeiaKhelghatibItems}, \ref{LankeitConsumptionItem} and \ref{MarrasViglialoroItem} explicitly provide conditions that involve $\chi \lVert v_0 \rVert_{L^\infty(\Omega)}$, exactly as in \eqref{Conditionv0chiTao} of \cite{TaoBoun}. Therefore, this term seems to be a very appropriate quantity that coherently characterizes the nature of models of the type \eqref{sysGeneralForConsumptionAndLogisticJohannes}. Consequently, it will play a role in our main result.
\end{remark}
\subsection{Gradient-dependent sources \texorpdfstring{\( h = h(u, \nabla u) \)}{h}}
Very recently, there has been growing interest in investigations focused on the dynamics and behavior of chemotactic systems that incorporate gradient-dependent sources. These studies seek to understand how external sources, influenced by the gradient of chemical concentrations, affect the movement and interaction of cells. We mention \cite{IshidaLankeitVigliloro-Gradientabbr,ACVgradient}, where real biological interpretations are also connected to ecological mechanisms. 

In particular, the source \( h = h(u, \nabla u) = \lambda u^\rho - \mu u^k - c |\nabla u|^\gamma \) describes a reaction with double dissipative action: one associated with the 0-order term (i.e., \( u^k \)) and the other with the 1-order term (i.e., \( |\nabla u|^\gamma \)). When introduced into chemotaxis models with production, this double damping action prevents blow-up situations; see, for example, \cite{IshidaLankeitVigliloro-Gradientabbr, ACVgradient} (and also \cite{VigliaGradTermDiffIntEqua}). 

In contrast, in the context of absorption chemotaxis systems, these terms allow for the boundedness of solutions even for small values of \( \mu \) and/or for large values of \( \lVert v_0 \rVert_{L^\infty(\Omega)} \). In this sense, since we are interested in its generalization, we mention the analysis of the following zero-flux consumption model with gradient-dependent source (as far as we know this is the only result so far available in the literature):
\begin{equation}\label{consumptionLinearWithGradient}
 \begin{cases}
 u_t = \Delta u - \chi \nabla \cdot u \nabla v + \lambda u - \mu u^2 - c | \nabla u |^\gamma, & \text{in } \Omega \times (0, \TM), \\
 v_t = \Delta v - uv, & \text{in } \Omega \times (0, \TM).
 \end{cases}
\end{equation}
In \cite{Columbu_arXiv_2024}, it is shown that the related initial-boundary value problem has a unique and uniformly bounded classical solution for \( \gamma \in \left( \frac{2n}{n+1}, 2 \right] \), regardless of the size of $\mu$ with respect $\chi \lVert v_0\rVert_{L^\infty(\Omega)}$. (Recall Remark \ref{RemarckChiV_0}.) 
\section{Main result and structure of the paper}\label{StructureMainReusltSection}
In agreement with all of the above, this paper is dedicated to the following problem
\begin{equation}\label{consum}
    \begin{cases}
        u_t = \nabla \cdot \tonda*{ (u+1)^{m_1-1}\nabla u  - \chi u(u+1)^{m_2-1}\nabla v} \\
       \hspace{7mm} + \lambda u - \mu u^2 - c \lvert \nabla u \rvert^\gamma & \text{in $\Omega \times (0,\TM)$},
       \\
        v_t = \Delta v - uv & \text{in $\Omega \times (0,\TM)$}, \\
        u_\nu = v_\nu = 0 & \text{on $\partial \Omega \times (0,\TM) $}, \\
        u(x,0) = u_0(x), \; v(x,0) = v_0(x) & x \in \bar\Omega,
    \end{cases}
\end{equation}
with $\Omega$ a bounded and smooth domain of $\R^n$ ($n\geq 2$), $m_1,m_2 \in \R$, $\lambda,\mu,c,\gamma$ proper positive reals and with  $\tmax\in(0,\infty]$.

As anticipated in some previous comments, the above model generalizes those in \eqref{consumptionLinearWithGradient} and in 
\eqref{sysGeneralForConsumptionAndLogisticJohannes}; in particular system \eqref{consumptionLinearWithGradient} is obtained by imposing $m_1=m_2=1$ in \eqref{consum} and the system studied in item \ref{MarrasViglialoroItem} is recovered for $c=0$.
We will make some considerations in this regard after having claimed our main result.
\subsection{Claim of the main result and structure of the paper}
In reference to our model under consideration (i.e, system \eqref{consum}), in order to present our result, we need to fix 
\begin{equation}\label{reginitialconditions}
\begin{cases}
\Omega\subset\R^n (n\geq 2) \;\textrm{a bounded domain of class}\;  C^{2+\delta},  \delta\in(0,1),\\
    u_0,v_0\colon \bar{\Omega}  \rightarrow \R^+ \mid u_0,  v_0 \in C_\nu^{2+\delta}(\bar\Omega)=\{\psi \in C^{2+\delta}(\bar{\Omega}): \psi_\nu=0 \textrm{ on }\partial \Omega\}.
    \end{cases}
\end{equation}
Having established the necessary preparations, we are now in a position to present the central result of this paper.
\begin{theorem}\label{teo}
 Let the hypotheses in \eqref{reginitialconditions} be fulfilled, $\lambda,\mu,c,\chi>0$,  $m_1,m_2 \in \R$. Then provided either 
   \begin{enumerate}[label={($\mathcal{A}_{\arabic*}$)}]
        \item \label{a1} $\displaystyle \max{\graffa*{\frac{2n}{n+1},\frac{n}{n+1}(2m_2-m_1+1)}}<\gamma\leq 2$
    \end{enumerate}
    or
    \begin{enumerate}[resume, label={($\mathcal{A}_{\arabic*}$)}]
        \item \label{a2} $\displaystyle \frac{n}{n+1}(2m_2-m_1+1)<\gamma\leq \frac{2n}{n+1}$ \qquad and \qquad $\displaystyle \mu>\frac{4 }{p_0}\mathcal{K} \norm{\chi v_0}_{L^{\infty}(\Omega)}^{2p_0}$,
    \end{enumerate}
    for some existing positive constants $p_0=p_0(m_1,m_2,n,\gamma)$ and $\mathcal{K}=\mathcal{K}(p_0)$, problem \eqref{consum} admits a unique and uniformly bounded classical solution, i.e., 
    \begin{equation*}
(u,v)\in \left(C^{2,1}( \Bar{\Omega} \times [0, \infty))\cap L^\infty(\Omega\times (0,\infty))\right)^2.
\end{equation*}
\end{theorem}
\begin{remark}
Let us make these comments.
\begin{enumerate}
    \item If in problem \eqref{consum} we set $m_1=m_2=1$, relations \ref{a1} and \ref{a2} recover what was established in \cite{Columbu_arXiv_2024}.
    \item The presence of the gradient nonlinearity in problem \eqref{consum}, i.e. $c>0$, allows one to obtain boundedness even for strong attraction effect (i.e., $m_2$ large), provided \ref{a1} is satisfied. This situation is different from the case where $c=0$; in this case $m_2<\frac{m_1+1}{2}$ and \( \mu \) large enough are required (see \cite{MarrasViglialoroMathNach}). 
\end{enumerate}
\end{remark}
The remainder of the paper is organized as follows. In $\S$\ref{PreliminariesSection}, we provide some preparatory and well-known preliminaries. Section $\S$\ref{ExistenceSolutionRegularizingSection} focuses on deriving a result concerning the local-in-time existence of classical solutions $(u,v)$ to system \eqref{consum} and some crucial properties related to the $u$- and $v$-components. Very importantly, in the same section, we will give a general  boundedness result for parabolic equations with gradient nonlinearities. 

Then, in $\S$\ref{AprioriEstimateSection}, we will control by means of a priori estimates the energy functional  
\[
\varphi(t) := \int_\Omega (u+1)^p + \chi^{2p}\int_\Omega |\nabla v|^{2p}\quad \textrm{on} \quad (0,\TM),
\]
defined, for $p>1$, in terms of the local solution. Subsequently, using a comparison principle, we provide a time-independent bound for $\varphi(t)$, which, in turn, gives bounds for $u$ in $L^p(\Omega)$ and $\nabla v$ in $L^{2p}(\Omega)$. Finally, we deduce that the local classical solution is actually bounded.
\section{Preliminaries}\label{PreliminariesSection}
The following are the technical and general results that will be referenced throughout the paper.
\begin{lemma} Let $p\in (1,\infty)$. Then for any $\psi\in C^2(\Bar{\Omega})$ satisfying $\psi\psi_{\nu}=0$ on $\partial \Omega$, the following inequality holds:
    \begin{equation}\label{grad4qquadro}
        \norm{\nabla \psi}^{2p+2}_{L^{2p+2}(\Omega)}\leq 2(4p^2+n)\norm{\psi}^2_{L^\infty(\Omega)} \norm*{\abs{\nabla \psi}^{p-1}D^2 \psi}^2_{L^2(\Omega)}.
    \end{equation}
    \begin{proof}
See  \cite[Lemma 2.2]{LankeitWangConsumptLogisticabbr}.
    \end{proof}
\end{lemma}
\begin{lemma}\label{lemmap} 
Let $n\in\N$ and 
\begin{equation*}
    \gamma>\max{\graffa*{\frac{2n}{n+1},\frac{n}{n+1}(2m_2-m_1+1)}},
\end{equation*}
 be satisfied. Then, there exists $p_1>1$ such that for every $p>p_1$ and
    \begin{table}[H]
    \centering
        \begin{subtable}[h]{0.66\textwidth}
        \begin{equation*}
             \bar{\theta}(p,\gamma)=\frac{\frac{p+\gamma-1}{\gamma}\tonda*{1-\frac{1}{p+1}}}{\frac{p+\gamma-1}{\gamma}+\frac 1n-\frac 1\gamma},
        \end{equation*}
        \end{subtable}
        \hfill
        \begin{subtable}[h]{0.33\textwidth}
        \begin{equation*}
            \bar{\sigma}(p,\gamma)=\frac{\gamma(p+1)}{p+\gamma-1},
        \end{equation*}
        \end{subtable}
        \\
        \begin{subtable}[h]{0.66\textwidth}
        \begin{equation*}
            \hat{\theta}(p)=\frac{\frac{p+\gamma-1}{\gamma}\tonda*{1-\frac{p}{(p+1)(p+2m_2-m_1-1)}}}{\frac{p+\gamma-1}{\gamma}+\frac 1n-\frac 1\gamma},
        \end{equation*}
        \end{subtable}
        \hfill
        \begin{subtable}[h]{0.33\textwidth}
        \begin{equation*}
            \tilde{\theta}(p)=\frac{\frac p2-\frac 12}{\frac p2-\frac 12 +\frac 1n},
        \end{equation*}             
        \end{subtable}
        \\
        \begin{subtable}[h]{0.66\textwidth}
        \begin{equation*}
          \hat{\sigma}(p)=\frac{\gamma(p+1)(p+2m_2-m_1-1)}{p(p+\gamma-1)}.  
        \end{equation*}             
        \end{subtable}
    \end{table}
these inequalities are complied: 
\begin{table}[H]
\centering
 \begin{subequations}
\begin{subtable}[h]{0.32\textwidth}
\begin{equation}\label{thetabar}
        0<\Bar{\theta}<1,
        \end{equation}
\end{subtable}
\hfill
\begin{subtable}[h]{0.32\textwidth}
\begin{equation}\label{thetasigmabar}
        \frac{\bar{\sigma}\bar{\theta}}{\gamma}<1.
        \end{equation}
\end{subtable}
\hfill
\begin{subtable}[h]{0.32\textwidth}
 \begin{equation}\label{thetahat} 
        0<\hat{\theta}<1,
    \end{equation}
\end{subtable}
\\
\begin{subtable}[h]{0.32\textwidth}
\begin{equation}\label{thetasigmahat}
        \frac{\hat{\sigma}\hat{\theta}}{\gamma}<1.
        \end{equation}
\end{subtable}
\hfill
\begin{subtable}[h]{0.32\textwidth}
\begin{equation}\label{thetatilde}
        0<\tilde{\theta}<1.
    \end{equation}
\end{subtable}
\end{subequations}
\end{table}
\begin{proof} 
The proof of the inequalities is straightforward as long as $p$ is chosen sufficiently large; let us indicate with $p_1$ such a value. We emphasize that the conditions $\gamma>\frac{2n}{n+1}$ and $\gamma>\frac{n}{n+1}(2m_2-m_1+1)$  are only needed to prove \eqref{thetasigmabar} and \eqref{thetasigmahat}, respectively.
\end{proof}
\end{lemma}
\section{Existence of local-in time solutions,  their properties and boundedness criterion}\label{ExistenceSolutionRegularizingSection}
The previous step toward establishing the boundedness of solutions to the system \eqref{consum} involves deriving a result that ensures the local in time existence of classical solutions. In particular, we aim to prove that solutions to the system exist within a certain time interval and are sufficiently smooth (i.e., they are classical solutions). In addition to this existence result, we outline some of the fundamental properties of these solutions. 
\begin{lemma}\label{theoremExistence}
 Let the hypotheses in \eqref{reginitialconditions} be fulfilled, $\lambda,\mu,c,\chi>0$,  $m_1,m_2 \in \R$ and $\gamma \geq 1$. Then 
there exist $\TM\in (0,\infty]$ and a unique couple of functions $(u,v$), with 
\begin{equation*}
(u,v)\in \left(C^{2+\delta,1+\frac{\delta}{2}}( \Bar{\Omega} \times [0, \TM))\right)^2
\end{equation*}
solving problem \eqref{consum}, and such that  \begin{equation}\label{dictomyCriteC2+del} 
 \text{if} \quad \TM<\infty \quad \text{then} \quad \lim_{t \to \TM} \left(\|u(\cdot,t)\|_{C^{2+\delta}(\bar\Omega)}+\|v(\cdot,t)\|_{C^{2+\delta}(\bar\Omega)}\right)=\infty.
 \end{equation}
 Additionally, we have
\begin{equation}\label{MaximumPricnipleRelation} 
u\geq 0\quad 0\leq v\leq \lVert v_0\rVert_{L^\infty(\Omega)}\quad \textrm{in}\quad \Omega \times (0,T_{max}), 
\end{equation}
\begin{equation}\label{boundednessmass}
        \into u \leq  \max{\graffa*{\into u_0(x)dx, \frac \lambda\mu \abs{\Omega}}} \quad \text{on $(0,\tmax)$},
    \end{equation}
and
\begin{equation}\label{boundednessgradv2}
    \nabla v\in L^\infty((0,\TM);L^2(\Omega)).
\end{equation}
\begin{proof}
Part of the proof is obtained by well-established methods
involving standard parabolic regularity theory in the framework of appropriate fixed point procedures, and the maximum principle; in the specific, see \cite{IshidaLankeitVigliloro-Gradientabbr,ACVgradient} (and also \cite{AmannBook} for alike results in more general settings). 

Further, by integrating the first equation and using the Neumann boundary conditions, we have
\begin{equation}\label{mass1}
    \frac{d}{dt}\into u\leq \lambda \into u -\mu \into u^2 \quad \text{for all $t\in(0,\tmax)$}.
\end{equation}
Moreover, an application of Hölder's inequality leads to 
\begin{equation*}
    \into u \leq \abs{\Omega}^{\frac 12}\tonda*{\into u^2}^{\frac 12} \quad \text{and, consequently, to} \quad -\into u^2\leq -\frac{1}{\abs{\Omega}}\tonda*{\into u}^2,
\end{equation*}
for every $t\in(0,\tmax)$.
By inserting the last estimate into \eqref{mass1} we obtain this initial problem
\begin{equation*}
\begin{dcases*}
    \frac{d}{dt}\int_\Omega u\leq \lambda \int_\Omega u - \frac{\mu}{\abs{\Omega}}\left(\int_\Omega u\right)^2 & \text{on $(0,\tmax)$},\\
    \left(\int_\Omega u \right)_{\big\vert_{t=0}}=\into u_0(x)dx,
\end{dcases*}
\end{equation*}
which ensures that $\int_\Omega u\leq \max{\graffa*{\into u_0(x)dx, \frac \lambda\mu \abs{\Omega}}}$ for all $t\in(0,\tmax)$.

Finally, relation \eqref{boundednessgradv2} is derived in  \cite[Lemma 3.4]{LankeitWangConsumptLogisticabbr}, with the support of \eqref{boundednessmass} itself.
\end{proof}
\end{lemma}
\subsection{Guaranteeing boundedness through estimates in Lebesgue spaces}
From now on, by $(u, v)$ we will refer to the unique local solution to model \eqref{consum}, which is defined in the domain $\Omega \times (0, T_{\text{max}})$ and provided by Lemma \ref{theoremExistence}. In the context of this lemma, the relation \eqref{dictomyCriteC2+del} holds, and in this specific case we say that the solution blows up at time $T_{\text{max}}$ in the $C^{2+\delta}(\Omega)$-norm. More precisely, it is important to note that this blow-up does not necessarily imply that the solution will blow up in the $L^\infty(\Omega)$-norm; in other words, while the solution may lose regularity and fail to stay smooth in the $C^{2+\delta}(\Omega)$-norm, it may still remain bounded in the $L^\infty(\Omega)$-norm.

To reformulate the extensibility criterion presented in relation \eqref{dictomyCriteC2+del}  in terms of uniform boundedness in time of the function $u$, we need to impose further restrictions on the class of gradient nonlinearities. Specifically, we require nonlinearity to satisfy conditions where the gradient term $\nabla u$ grows at most quadratically, that is, the nonlinearity is of the form $|\nabla u|^2$ at most. This restriction is precisely what allows us to establish a well-posedness framework in which we can discuss the uniform-in-time boundedness of $u$ and avoid the blow-up in $L^\infty(\Omega)$-norm.
\begin{lemma}[Boundedness criterion]\label{ExtCriterionLemma} Let $\gamma \in [1,2]$. There is $p_2>1$ such that  if $u \in L^\infty((0, \TM); L^{p_2}(\Omega))$ then $
u \in L^\infty((0, \infty); L^\infty(\Omega))$.
\begin{proof}
First, we observe that for any \( \lambda, \mu, c > 0 \) and  \( \gamma \geq 1 \), the function \( u \mapsto \lambda u - \mu u^2 - c |\nabla u|^\gamma \), for \( u \geq 0 \), is bounded from above by a positive constant \( L \); henceforth, we can write
\[
g(x, t) := \lambda u- \mu u^2 - c |\nabla u|^\gamma \leq L \quad \text{in} \quad \Omega \times (0, T_{\text{max}}).
\]
Subsequently, with a view to \cite[Appendix A]{TaoWinkParaPara}, and consistently with that symbology, it is seen that if \( u \) classically solves in \( \Omega \times (0, T_{\text{max}}) \) the first equation in problem \eqref{consum}, it also solves \cite[problem (A.1)]{TaoWinkParaPara} with
\[
D(x, t, u) = (u + 1)^{m_1 - 1}, \quad f(x, t) = -\chi u(u + 1)^{m_2 - 1} \nabla v, \quad g(x, t) = L.
\]
Thereafter, by choosing an appropriately large \( p_2>1 \) such that from the hypothesis $u \in L^\infty((0, \TM); L^{p_2}(\Omega))$ parabolic regularity results applied to the equation \( v_t = \Delta v - uv \) give the necessary regularity for \( \nabla v \), we have that (A2)--(A10) are all satisfied, and henceforth we have by applying \cite[Lemma A.1]{TaoWinkParaPara} that actually $u \in L^\infty((0, \TM); L^{\infty}(\Omega))$.

For the details of the implication $$u \in L^\infty((0, \TM); L^{\infty}(\Omega)) \Rightarrow  u \in L^\infty((0, \infty); L^{\infty}(\Omega))$$ 
 we refer to  \cite[Lemma 5.2]{ACVgradient}; herein we just mention that by exploiting in the second equation of model \eqref{consum} the boundedness of $u$ in $L^{\infty}(\Omega)$ for all $t\in [0,\TM]$, we have that of \( \| v(\cdot, t) \|_{C^{1+\delta}(\Omega)} \), and in turn thanks to $\gamma\in[1,2]$, the uniform-in-time boundedness of $\norm*{u(\cdot,t)}_{C^{2+\delta}(\bar{\Omega})}+\norm*{v(\cdot,t)}_{C^{2+\delta}(\bar{\Omega})}$; so the conclusion follows by contradiction thanks to the extensibility criterion \eqref{dictomyCriteC2+del}.
\end{proof}
\end{lemma}
\section{Uniform bounds and proof of Theorem \ref{teo}}\label{AprioriEstimateSection}
To effectively apply Lemma \ref{ExtCriterionLemma}, we first need to derive some a priori time-independent estimates. These estimates will provide crucial bounds for $u$ and $v$, essential for the correct application of the lemma mentioned. To this end, we will study for any $p>1$ the time evolution of the functional
\begin{equation}\label{DefiFunctional}
\varphi(t) := \int_{\Omega} (u + 1)^p  + \chi^{2p}\int_{\Omega} |\nabla v|^{2p} \quad \textrm{for all} \quad t\in (0,\TM),
\end{equation}
the goal being to establish its boundedness. This will be achieved by deriving a differential problem of the type
\begin{equation}\label{AbsortpionProblem}
\begin{cases}
    \varphi'(t)+\varphi(t)\leq C \quad \textrm{for all} \; t\in (0,\TM),\\
    \varphi(0)=\int_{\Omega} (u_0(x) + 1)^p dx  + \chi^{2p}\int_{\Omega} |\nabla v_0(x)|^{2p}dx,
\end{cases}
\end{equation}
for some $C>0$, which by an ODI comparison reason entails 
\[\int_\Omega (u+1)^p+\chi^{2p}\int_\Omega |\nabla v|^{2p}\leq \max\{C,\varphi(0)\} \quad \textrm{for all } t \in (0,\TM).\]
\textit{Based on the points discussed up to this point, in the following sections, our objective is to establish an inequality for the function \( \varphi \) as that  presented in  \eqref{AbsortpionProblem}.}
\subsection{A priori estimates}
 The next computations collect uniform-in-time bounds of some norms of $u$ and $v$;  some related derivations will involve the continuous function
    \begin{equation}\label{kappa2}
        K(p,n,\eta)\coloneqq\frac{2p^{\frac{p+1}{2}}(p+n+\eta-1)^{\frac{p+1}{2}}}{p+1}\tonda*{\frac{8(4p^2+n)(p-1)}{p(p+1)}}^{\frac{p-1}{2}},
    \end{equation}
defined for all $p>1$ and $\eta>0$,
as well as
\begin{equation}\label{cUnder-cSuper} 
  K_1=p(p + n - 1 + \eta) \|v_0\|^2_{L^\infty(\Omega)}\quad \textrm{and}\quad K_2=K(p,n,\eta)\nivz^{2p}.  
\end{equation}
All the constant $c_i$, $i=1,2,3,\ldots$ appearing below are positive.
\begin{lemma}
Let the hypotheses of Lemma \ref{lemmap} be fulfilled, and $p_1>1$ the value therein obtained.  Then for every $p>p_1$ and for all $\varepsilon,\hat{c}>0$ these inequalities hold:
    \begin{equation}\label{gammamaggiore}
         \hat{c}\into(u+1)^{p+1}\leq \varepsilon \into\abs*{\nabla (u+1)^\frac{p+\gamma-1}{\gamma}}^\gamma + \const{sg}  \quad \text{on $(0,\tmax)$},
    \end{equation}
     \begin{equation}\label{gammamaggiore2}
         \hat{c}\into(u+1)^{\frac{(p+1)(p+2m_2-m_1-1)}{p}}\leq \varepsilon \into\abs*{\nabla (u+1)^\frac{p+\gamma-1}{\gamma}}^\gamma + \const{sg1}  \quad \text{on $(0,\tmax)$},
    \end{equation}
    and
    \begin{equation}\label{gngradv2}
        \into \abs{\nabla v}^{2p}\leq \frac p4 \into \abs{\nabla v}^{2p-2}\abs*{D^2 v}^2+\const{b7} \quad \text{on $(0,\tmax)$}.
    \end{equation}
\begin{proof}
In order to prove relation \eqref{gammamaggiore}, we use the Gagliardo--Nirenberg inequality (see \cite{Nirenber_GagNir_Ineque}) in conjunction with \eqref{thetabar} and \eqref{boundednessmass}, obtaining (we set for commodity in the next few lines $\hat{u}=(u+1)$)
\begin{equation*}
\begin{split}
   \hat{c} \into (u+1)^{p+1}=&  \norm*{\hat{u}^\frac{p+\gamma-1}{\gamma}}_{L^{\frac{\gamma(p+1)}{p+\gamma-1}}(\Omega)}^\frac{\gamma(p+1)}{p+\gamma-1}\\
    \leq& \const{s2}\tonda*{\norm*{\nabla\hat{u}^\frac{p+\gamma-1}{\gamma}}_{L^\gamma(\Omega)}^{\Bar{\theta}}  \norm*{\hat{u}^\frac{p+\gamma-1}{\gamma}}_{L^\frac{\gamma}{p+\gamma-1}(\Omega)}^{1-\Bar{\theta}}+\norm*{\hat{u}^\frac{p+\gamma-1}{\gamma}}_{L^\frac{\gamma}{p+\gamma-1}(\Omega)}}^{\Bar{\sigma}}\\
    \leq& \const{s3}\tonda*{\into\abs*{\nabla \hat{u}^\frac{p+\gamma-1}{\gamma}}^\gamma}^\frac{\Bar{\sigma}\Bar{\theta}}{\gamma}+\const{s4} \leq \varepsilon \into\abs*{\nabla \hat{u}^\frac{p+\gamma-1}{\gamma}}^\gamma + \const{sg}\quad \text{on }(0,\tmax),
\end{split}
\end{equation*}
where in the last estimate we applied Young's inequality, justified by relation \eqref{thetasigmabar}, and $(A+B)^s\leq \max\{1, 2^{s-1}\}(A^s+B^s)$, valid for any $A,B\geq 0$ and $s>0$. (We will tacitly employ such an algebraic inequality without mentioning it.)

As to relation \eqref{gammamaggiore2}, relying on \eqref{thetahat} and \eqref{thetasigmahat}, similar rearrangements give on $(0,\tmax)$
\begin{equation*}
\begin{split}
    \hat{c} \into (u+1)^{\frac{(p+1)(p+2m_2-m_1-1)}{p}}=&\norm*{\hat{u}^\frac{p+\gamma-1}{\gamma}}_{L^{\frac{\gamma(p+1)(p+2m_2-m_1-1)}{p(p+\gamma-1)}}(\Omega)}^\frac{\gamma(p+1)(p+2m_2-m_1-1)}{p(p+\gamma-1)}\\
    \leq& \const{s2}\left( \norm*{\nabla\hat{u}^\frac{p+\gamma-1}{\gamma}}_{L^\gamma(\Omega)}^{\hat{\theta}}  \norm*{\hat{u}^\frac{p+\gamma-1}{\gamma}}_{L^\frac{\gamma}{p+\gamma-1}(\Omega)}^{1-\hat{\theta}} \right. \\ & \left. \hspace{40mm}+\norm*{\hat{u}^\frac{p+\gamma-1}{\gamma}}_{L^\frac{\gamma}{p+\gamma-1}(\Omega)}\right)^{\hat{\sigma}}\\
    \leq& \const{s3}\tonda*{\into\abs*{\nabla \hat{u}^\frac{p+\gamma-1}{\gamma}}^\gamma}^\frac{\hat{\sigma}\hat{\theta}}{\gamma}+\const{s4} \leq \varepsilon \into\abs*{\nabla \hat{u}^\frac{p+\gamma-1}{\gamma}}^\gamma + \const{sg1}.
\end{split}
\end{equation*}
Ultimately, by recalling \eqref{boundednessgradv2} and \eqref{thetatilde} we have
    \begin{equation*}
    \begin{split}
        \into \abs{\nabla v}^{2p}=&\into \tonda*{\abs{\nabla v}^p}^2=\norm*{\abs{\nabla v}^p}^2_{L^2(\Omega)}\\
        \leq & \const{b10} \tonda*{\norm*{\nabla \abs*{\nabla v}^p}^{\tilde{\theta}}_{L^2(\Omega)}\norm*{\abs{\nabla v}^p}^{1-\tilde{\theta}}_{L^{\frac 2p}(\Omega)}+\norm*{\abs*{\nabla v}^p}_{L^{\frac 2p}(\Omega)}}^2\\
        \leq & \frac{1}{4p} \into \abs*{\nabla \abs*{\nabla v}^p}^2 +\const{b11} \leq \frac p4 \into \abs*{\nabla v}^{2p-2}\abs*{D^2 v}^2 +\const{b7} \quad \text{on $(0,\tmax)$},
    \end{split}
    \end{equation*}
so that also relation \eqref{gngradv2} is achieved.
\end{proof}
\end{lemma}
\begin{lemma}
    For every $p>1$ and for all $\hat{c}>0$ we have
    \begin{equation}\label{youngv0}
        K_1 \into u^2\absnv^{2p-2}\leq \frac p4 \into \absnv^{2p-2}\dqvq+K_2\into u^{p+1} \quad \text{on $(0,\tmax)$},
    \end{equation}
    \begin{equation}\label{up1grav2}
    \begin{split}
        \hat{c}\into (u+1)^{p+2m_2-m_1-1}\abs{\nabla v}^2\leq & \frac p4 \chi^{2p} \into \abs{\nabla v}^{2p-2}\abs*{D^2 v}^2 \\ & + \const{gh} \into (u+1)^{\frac{(p+1)(p+m_2-m_1-1)}{p}} \quad \text{on $(0,\tmax)$},
    \end{split}
    \end{equation}
where  $K_1$ and $K_2$ have been introduced in \eqref{cUnder-cSuper}.
\begin{proof}
Let us prove  bounds \eqref{youngv0} and \eqref{up1grav2}. Indeed, an application of Young's inequality, in conjunction with \eqref{grad4qquadro} and \eqref{MaximumPricnipleRelation}, leads on $(0,\tmax)$ and for any $\varepsilon_1>0$ and $C(\varepsilon_1)=\frac{2 K_1}{p+1}\left(\varepsilon_1\frac{p+1}{K_1(p-1)}\right)^\frac{1-p}{2}$ to
\begin{equation*}
\begin{split}
     K_1 \into u^2\absnv^{2p-2}\leq & \varepsilon_1 \into \absnv^{2(p+1)}+C(\varepsilon_1)\into u^{p+1}\\
    \leq &\varepsilon_1 2(4p^2+n)\nivz^2 \into \absnv^{2p-2}\dqvq \\
     &+ \frac{2}{p+1}(\varepsilon_1\frac{p+1}{p-1})^{-\frac{p-1}{2}}K_1^{\frac{p+1}{2}}\into u^{p+1} \\
    = & \frac p4 \into \absnv^{2p-2}\dqvq+K_2\into u^{p+1},
\end{split}
\end{equation*}
and for $\varepsilon_1=\frac{p}{8(4p^2+n)\nivz^2}$ relation \eqref{youngv0} is attained. Likewise, we have for all $t\in(0,\tmax)$,
\begin{equation*}
    \begin{split}
    \hat{c}\into (u+1)^{p+2m_2-m_1-1}\abs{\nabla v}^2\leq 
    & \frac p4 \chi^{2p}\into \abs{\nabla v}^{2p-2}\abs*{D^2 v}^2 \\ & + \const{gh} \into (u+1)^{\frac{(p+1)(p+m_2-m_1-1)}{p}},
    \end{split}
\end{equation*}
and \eqref{up1grav2} is established as well.
\end{proof}
\end{lemma}
\begin{lemma}
For any $p > 1$ and for every $t\in(0,\tmax)$ we have
\begin{equation}\label{up}
\begin{split}
    \frac{d}{dt}\into (u+1)^p\leq & \const{fa1}\into (u+1)^{p+2m_2-m_1-1}\abs{\nabla v}^2+\const{fa2}\into (u+1)^p\\
    &-\frac{\mu p}{4}\into (u+1)^{p+1}-\const{fa4}\into \abs*{\nabla u^{\frac{p+\gamma-1}{\gamma}}}^\gamma+\const{fa5}.
\end{split}
\end{equation}
\begin{proof}
We compute the time derivative of $\int_\Omega  (u+1)^p$ and apply Young's inequality, yielding on $(0, \TM)$
\begin{equation}\label{IneqA1}
\begin{split}
        \frac{d}{dt}\into (u+1)^p=&-p(p-1)\into (u+1)^{p+m_1-3}\abs{\nabla u}^2\\
        &+\chi p(p-1)\into u(u+1)^{p+m_2-3}\nabla u\cdot \nabla v +\lambda p \into (u+1)^p \\
        & -\mu p \into u^2(u+1)^{p-1} -c p \into (u+1)^{p-1} \abs{\nabla u}^\gamma.
\end{split}
\end{equation}
In turn from the inequality $(u+1)^2\leq 2(u^2+1)$ we have through  Young's inequality that on $(0,\TM)$ one has
\begin{equation*}
\begin{split}
    -\mu p \into u^2(u+1)^{p-1}\leq & -\frac{\mu p}{4}\into (u+1)^{p+1}+\mu p\into (u+1)^{p-1}\\
    \leq &
    -\frac{\mu p}{4}\into (u+1)^{p+1}+\const{fg1}\into (u+1)^p +\const{fa5},
\end{split}
\end{equation*}
and also, by the Young inequality again,
\begin{equation*}
\begin{split}
    \chi p (p-1)\into (u+1)^{p+m_2-2}\nabla u \cdot \nabla v \leq &\frac{\const{fb1}}{2}\into (u+1)^{p+m_1-3}\abs{\nabla u}^2 \\
    &+\const{fa1}\into (u+1)^{p+2m_2-m_1-1}|\nabla v|^2.
\end{split}
\end{equation*}
These last estimates inserted in \eqref{IneqA1} give the claim.
\end{proof}
\end{lemma}
\begin{lemma}
For all $\eta>0$,  $p\in(1,\infty)$ and  $t\in(0,\tmax)$ this inequality is satisfied
\begin{equation}\label{v2p}
    \frac{d}{dt}\into \absnv^{2p}+p\into \absnv^{2p-2} \dqvq \leq K_1 \into u^2\absnv^{2p-2}+\const{c1},
\end{equation}
with $K_1$ given in \eqref{cUnder-cSuper}.
\begin{proof}
The proof is contained in \cite[Lemma 4.2]{LankeitWangConsumptLogisticabbr}.
\end{proof}
\end{lemma}
\begin{lemma}\label{lemmanormap} Let either 
\begin{equation}\label{ConditionA1}
    \gamma>\max{\graffa*{\frac{2n}{n+1},\frac{n}{n+1}(2m_2-m_1+1)}},
\end{equation}
or, for any $\eta>0$ and all $p>p_1$ ($p_1$ being the value detected in Lemma \ref{lemmap}) and $K(p,n,\eta)$ as in \eqref{kappa2},
\begin{equation}\label{condsommacmuconetaep}
    \gamma>\frac{n}{n+1}(2m_2-m_1+1) \quad \text{and} \quad \mu>\frac 4p K(p,n,\eta)\norm{\chi v_0}_{L^{\infty}(\Omega)}^{2p}.
\end{equation}
Then  $u \in L^\infty((0,\TM);L^p(\Omega))$, for all $p>p_1$.
\begin{proof}
Let $\varphi(t)$ fixed as in \eqref{DefiFunctional} and recall the definition of $K_1$ and $K_2$ in \eqref{cUnder-cSuper};  by exploiting bounds  \eqref{up} and \eqref{v2p}, we have that for all $t\in (0,\TM)$ one has 
    \begin{equation*}
    \begin{split}
        \varphi'(t)+\varphi(t)\leq& \const{fa1}\into (u+1)^{p+2m_2-m_1-1}\abs{\nabla v}^2+\const{fa23}\into (u+1)^p-\frac{\mu p}{4}\into (u+1)^{p+1} \\
        & -\const{fa4}\into \abs*{\nabla u^{\frac{p+\gamma-1}{\gamma}}}^\gamma -\chi^{2p}p\into \absnv^{2p-2} \dqvq \\
        &  + \chi^{2p}K_1 \into u^2\absnv^{2p-2}+\chi^{2p}\into \absnv^{2p}+\const{fa56},
    \end{split}
    \end{equation*}
   which for all $t \in (0,\TM)$ and any $\rho>0$ is simplified, by virtue of \eqref{gngradv2}, \eqref{youngv0},  \eqref{up1grav2} and the Young inequality, into 
   \begin{equation*}
    \begin{split}
    \varphi'(t)+\varphi(t)\leq &\tonda*{K_2\chi^{2p}-\frac{\mu p}{4}+\rho}\into (u+1)^{p+1} 
     \\ & +\const{gh}\into (u+1)^{\frac{(p+1)(p+m_2-m_1-1)}{p}}-\const{fa4}\into \abs*{\nabla u^{\frac{p+\gamma-1}{\gamma}}}^\gamma +\const{gh24}.
    \end{split}
\end{equation*}
In turn, if relation \eqref{ConditionA1} is satisfied, independently by the sing of $K_2\chi^{2p}-\frac{\mu p}{4}+\rho$,  we use \eqref{gammamaggiore} and  \eqref{gammamaggiore2}, with $\varepsilon=\frac{\const{fa4}}{2}$, so to obtain an absorption problem as \eqref{AbsortpionProblem}. Conversely, under assumption \eqref{condsommacmuconetaep}, we can choose $\rho$ sufficiently small so to have $K_2\chi^{2p}-\frac{\mu p}{4}+\rho\leq 0$;
subsequently, again \eqref{gammamaggiore2} leads to an alike problem as that in \eqref{AbsortpionProblem}. Definitively, we can conclude.
\end{proof}
\end{lemma}
\subsection{Proof of Theorem \ref{teo}}
We have gathered all the essential elements and mathematical tools required to  establish our result.
\begin{proof}
For $p_1$ and $p_2$ the constants from Lemmas \ref{lemmap} and \ref{ExtCriterionLemma}, respectively, let us set $p_0=\max\{p_1,p_2\}$. (Let us observe that $p_0=p_0(m_1,m_2,n,\gamma)$.)

If relation \ref{a1} is valid, Lemmata \ref{lemmanormap} and \ref{ExtCriterionLemma} directly ensure the claim. Conversely, let $\mathcal{K}(p_0)=K(p_0,n,0)$, being $K(p,n,\eta)$ the function in \eqref{kappa2} and let assumption \ref{a2} be satisfied. Since $K(p,n,\eta)$ is a continuous function of the variables $p$ and $\eta$, by continuity arguments  there exist $p>p_0$ and $\eta>0$ such that the second condition in \eqref{condsommacmuconetaep} is fulfilled; consequently, the same  mentioned lemmas give the conclusion.
\end{proof}
\subsubsection*{Acknowledgements}
DAS has been supported by UCA FPU contract UCA/REC14VPCT/2020 funded by Universidad de Cádiz, by mobility grants funded by Plan Propio - UCA 2022-2023 and by a Graduate Scholarship funded by The University of Tennessee at Chattanooga. AC and GV are members of the {\em Gruppo Nazionale per l’Analisi Matematica, la Probabilità e le loro Applicazioni} (GNAMPA) of the Istituto Nazionale di Alta Matematica (INdAM).
GV is also supported by MIUR (Italian Ministry of Education, University and Research) Prin 2022 \textit{Nonlinear differential problems with applications to real phenomena} (Grant Number: 2022ZXZTN2). AC is also supported by GNAMPA-INdAM Project \textit{Problemi non lineari di tipo stazionario ed evolutivo} (CUP--E53C23001670001).

\end{document}